# ESTIMATION OF FUNCTIONAL DERIVATIVES

By Peter Hall[1], Hans-Georg Müller[2] and Fang Yao[3]

*University of Melbourne, University of California, Davis
and University of Toronto*

Situations of a functional predictor paired with a scalar response are increasingly encountered in data analysis. Predictors are often appropriately modeled as square integrable smooth random functions. Imposing minimal assumptions on the nature of the functional relationship, we aim to estimate the directional derivatives and gradients of the response with respect to the predictor functions. In statistical applications and data analysis, functional derivatives provide a quantitative measure of the often intricate relationship between changes in predictor trajectories and those in scalar responses. This approach provides a natural extension of classical gradient fields in vector space and provides directions of steepest descent. We suggest a kernel-based method for the nonparametric estimation of functional derivatives that utilizes the decomposition of the random predictor functions into their eigenfunctions. These eigenfunctions define a canonical set of directions into which the gradient field is expanded. The proposed method is shown to lead to asymptotically consistent estimates of functional derivatives and is illustrated in an application to growth curves.

**1. Introduction.** Situations where one is given a functional predictor and a continuous scalar response are increasingly common in modern data analysis. While most studies to date have focused on functional linear models, the structural constraints imposed by these models are often undesirable. To enhance flexibility, several nonparametric functional regression approaches have been discussed. Since these models are not subject to any assumptions except smoothness, they are very widely applicable. The price one pays, of

Received June 2008; revised January 2009.
[1]Supported in part by an Australian Research Council Fellowship.
[2]Supported in part by NSF Grants DMS-05-05537 and DMS-08-06199.
[3]Supported in part by a NSERC Discovery grant.

course, is that convergence will be slower when compared with functional linear models. The situation is comparable to that of extending ordinary linear regression to nonparametric regression. By abandoning restrictive assumptions, such extensions greatly enhance flexibility and breadth of applicability. Under suitable regularity assumptions, convergence of such functional nonparametric models is guaranteed for a much larger class of functional relationships, and this insurance is often well worth the slower rates of convergence.

Suppose we observe a sample of i.i.d. data $(X_1, Y_1), \ldots, (X_n, Y_n)$, generated by the model

$$Y = g(X) + \varepsilon, \tag{1}$$

where $X$ is a random function in the class $L_2(\mathcal{I})$ of square-integrable functions on the interval $\mathcal{I} = [0, 1]$, $g$ is a smooth functional from $L_2(\mathcal{I})$ to the real line and $\varepsilon$ represents an error, independent of $X$, with zero expected value and finite variance. In the nonparametric approach, one aims to conduct inference about $g$ without imposing specific structure, usually that $g$ is a linear functional. The traditional functional linear model would have $g(x) = a + \int bx$, where $a$ is a constant and $b$ a function, but even here the "regression parameter function" $b$ cannot be estimated at the parametric rate $n^{-1/2}$, unless it is subject to a finite-parameter model; this model has been well investigated in the literature. Examples of such investigations include Ramsay and Dalzell (1991), Cuevas, Febrero and Fraiman (2002), Cardot et al. (2003), Cardot, Ferraty and Sarda (2003), James and Silverman(2005), Ramsay and Silverman (2005), Yao, Müller and Wang (2005b) and Hall and Horowitz (2007).

While the functional linear regression model has been shown to provide satisfactory fits in various applications, it imposes a linear restriction on the regression relationship and, therefore, cannot adequately reflect nonlinear relations. The situation is analogous to the case of a simple linear regression model where a nonparametric regression approach often provides a much more adequate and less biased alternative approach. Likewise, there is sometimes strong empirical evidence, for example, in the form of skewness of the distributions of empirical component scores, that the predictor function $X$ is not Gaussian. The problem of estimating a nonparametric functional regression relation $g$ in the general setting of (1) is more difficult compared to functional linear regression, and the literature is much sparser. It includes the works of Gasser, Hall and Presnell (1998) and Hall and Heckman (2002) on the estimation of distributions and modes in function spaces, and of Ferraty and Vieu (2003, 2004, 2006) on nonparametric regression with functional predictors. Recent developments are reviewed in Ferraty, Mas and Vieu (2007).

To lay the foundations for our study, we introduce an orthonormal basis for $L_2(\mathcal{I})$, say $\psi_1, \psi_2, \ldots$, which, in practice, would generally be the basis connected to the spectrum of the covariance operator, $V(s,t) = \text{cov}\{X(s), X(t)\}$:

$$V(s,t) = \sum_{j=1}^{\infty} \theta_j \psi_j(u) \psi_j(v), \tag{2}$$

where the $\psi_j$'s are the orthonormal eigenfunctions, and the $\theta_j$'s are the respective eigenvalues of the linear operator with kernel $V$. The terms in (2) are ordered as $\theta_1 \geq \theta_2 \geq \cdots$. The empirical versions of the $\psi_j$'s and $\theta_j$'s arise from a similar expansion of the standard empirical approximation $\widehat{V}$ to $V$,

$$\widehat{V}(s,t) = \frac{1}{n} \sum_{i=1}^{n} \{X_i(s) - \bar{X}(s)\}\{X_i(t) - \bar{X}(t)\} = \sum_{j=1}^{\infty} \hat{\theta}_j \hat{\psi}_j(s) \hat{\psi}_j(t), \tag{3}$$

where $\bar{X} = n^{-1} \sum_i X_i$ and order is now determined by $\hat{\theta}_1 \geq \hat{\theta}_2 \geq \cdots$. The eigenvalues $\hat{\theta}_j$ vanish for $j \geq n+1$, so the functions $\hat{\psi}_{n+1}, \hat{\psi}_{n+2}, \ldots$ may be determined arbitrarily.

The centered form of $X$ admits a Karhunen–Loève expansion

$$X - E(X) = \sum_{j=1}^{\infty} \xi_j \psi_j, \tag{4}$$

where the principal components $\xi_j = \int_{\mathcal{I}} X \psi_j$ are uncorrelated and have zero means and respective variances $\theta_j$. Their empirical counterparts are computed using $\hat{\psi}_j$ in place of $\psi_j$.

The paper is organized as follows. In Section 2, we describe the kernel-based estimators that we consider for estimating the nonparametric regression function $g$ in model (1) on the functional domain and for estimating functional derivatives in the directions of the eigenfunctions $\psi_j$. In Section 3, rates of convergence for kernel estimators $\hat{g}$ of the nonparametric regression function $g$ are obtained under certain regularity assumptions on predictor processes and their spectrum (Theorems 1 and 2). These results then lead to the consistency property (Theorem 3) for functional derivatives. A case study concerning an application of functional derivatives to the Berkeley longitudinal growth study is the theme of Section 4, followed by a compilation of the proofs and additional results in Section 5.

**2. Proposed estimation procedures.** Define the Nadaraya–Watson estimator

$$\hat{g}(x) = \frac{\sum_i Y_i K_i(x)}{\sum_i K_i(x)},$$

where $K_i(x) = K(\|x - X_i\|/h)$, $K$ is a kernel function and $h$ a bandwidth. Here, $\|\cdot\|$ denotes the standard $L^2$ norm. Similar kernel estimators have

been suggested in the literature. We refer to Ferraty and Vieu (2006) for an overview regarding these proposals and also for the previously published consistency results for the estimation of $g$. While the focus of this paper is on the estimation of functional derivatives in the general framework of model (1), using the spectral decomposition for predictor processes $X$ and characterizing these processes by their eigenbasis also leads to useful and relevant results regarding the estimation of $g$. These results are given in Theorems 1 and 2 below, while Theorem 4 provides relevant bounds for the probability that $X$ lies in a small ball, and Theorem 3 yields the desired asymptotic consistency of the proposed functional derivative estimator defined at (7).

For simplicity, we shall suppose the following (although more general conditions may be imposed).

ASSUMPTION 1. Kernel $K$ is nonincreasing on $[0, c]$, where $c > 0$, and the support of $K$ equals $[0, c]$.

The derivative of $g$ at $x$ is defined to be the linear operator $g_x$ with the property that, for functions $y$ and scalars $\delta$,

$$g(x + \delta y) = g(x) + \delta g_x y + o(\delta)$$

as $\delta \to 0$. We may write

$$(5) \qquad g_x = \sum_{j=1}^{\infty} \gamma_{xj} t_j,$$

where $\gamma_{xj} = g_x \psi_j$ is a scalar, and $t_j$ denotes the operator that takes $y$ to $y_j = t_j(y) = \int y \psi_j$. We can think of $\gamma_{xj}$ as the component of $g_x$ in the direction $\psi_j$.

From knowledge of the operator $g_x$, accessible through the components $\gamma_{xj}$, we can obtain information about functional gradients and extrema. For example, suppose $a_x^{\min} = (a_{x1}^{\min}, a_{x2}^{\min}, \ldots)$ and $a_x^{\max} = (a_{x1}^{\max}, a_{x2}^{\max}, \ldots)$ are defined as the vectors $a = (a_1, a_2, \ldots)$ that, respectively, minimize and maximize $|g_x a|$, where

$$(6) \qquad g_x a = \sum_{j=1}^{\infty} \gamma_{xj} a_j,$$

over functions $a = \sum_j a_j \psi_j$ for which $\|a\| = 1$ (i.e., such that $\sum_j a_j^2 = 1$). Then, the function $g$ changes fastest as we move away from $x$ in the direction of $a_x^{\max} = \sum_j a_{xj}^{\max} \psi_j$, which, therefore, is a gradient direction. The function changes least when we move from $x$ in the direction of $a_x^{\min} = \sum_j a_{xj}^{\min} \psi_j$. Extremal points are characterized by $\gamma_{xj} = 0$ for all $j$, and their identification is of obvious interest to identify predictor functions associated with maximal or minimal responses, and also the level of these responses.

Thus, the components $\gamma_{xj}$ are of intrinsic interest. As a prelude to estimating them, we introduce $Y_{i_1 i_2} = Y_{i_1} - Y_{i_2}$ and $\hat{\xi}_{i_1 i_2 j} = \int_\mathcal{I}(X_{i_1} - X_{i_2})\hat{\psi}_j$, the latter being an empirical approximation to $\xi_{i_1 i_2 j} = \xi_{i_1 j} - \xi_{i_2 j}$ (i.e., to the difference between the principal components $\xi_{ij} = \int X_i \psi_j$ for $i = i_1, i_2$). Define

$$Q_{i_1 i_2 j} = 1 - \frac{|\int (X_{i_1} - X_{i_2})\hat{\psi}_j|^2}{\|X_{i_1} - X_{i_2}\|^2} = 1 - \frac{\hat{\xi}_{i_1 i_2 j}^2}{\|X_{i_1} - X_{i_2}\|^2},$$

which represents the proportion of the function $X_{i_1} - X_{i_2}$, that is, "not aligned in the direction of $\hat{\psi}_j$." Therefore, $Q_{i_1 i_2 j}$ will be small in cases where $X_{i_1} - X_{i_2}$ is close to being in the direction of $\hat{\psi}_j$, and will be larger in other settings. We suggest taking

$$(7) \qquad \hat{\gamma}_{xj} = \frac{\sum \sum_{i_1, i_2}^{(j)} Y_{i_1 i_2} K(i_1, i_2, j | x)}{\sum \sum_{i_1, i_2}^{(j)} \hat{\xi}_{i_1 i_2 j} K(i_1, i_2, j | x)}.$$

Here, $\sum \sum_{i_1, i_2}^{(j)}$ denotes summation over pairs $(i_1, i_2)$ such that $\hat{\xi}_{i_1 i_2 j} > 0$,

$$(8) \qquad K(i_1, i_2, j | x) = K\left(\frac{\|x - X_{i_1}\|}{h_1}\right) K\left(\frac{\|x - X_{i_2}\|}{h_1}\right) K\left(\frac{Q_{i_1 i_2 j}}{h_2}\right),$$

$K$ is a kernel function and $h_1$ and $h_2$ denote bandwidths. On the right-hand side of (8), the last factor serves to confine the estimator's attention to pairs $(i_1, i_2)$, for which $X_{i_1} - X_{i_2}$ is close to being in the direction of $\hat{\psi}_j$, and the other two factors restrict the estimator to $i_1$ and $i_2$, such that both $X_{i_1}$ and $X_{i_2}$ are close to $x$. The estimator $\hat{\gamma}_{xj}$ uses two smoothing parameters, $h_1$ and $h_2$.

## 3. Theoretical properties.

3.1. *Consistency and convergence rates of estimators of $g$.* To ensure consistency, we ask that the functional $g$ be continuous at $x$ (i.e., that for functions $y$ and scalars $\delta$, the following holds).

ASSUMPTION 2.

$$(9) \qquad \sup_{y:\|y\|\leq 1} |g(x + \delta y) - g(x)| \to 0 \qquad \text{as } \delta \downarrow 0,$$

and the bandwidth $h$ does not decrease to zero too slowly, in the sense that, with $c$ as in Assumption 1,

$$(10) \qquad h = h(n) \to 0 \quad \text{and} \quad nP(\|X - x\| \leq c_1 h) \to \infty \qquad \text{as } n \to \infty,$$

where $c_1 = c$ if $K(c) > 0$, and otherwise $c_1 \in (0, c)$.

Given $C > 0$, $x \in L_2(\mathcal{I})$ and $\alpha \in (0,1]$, let $\mathcal{G}(C, x, \alpha)$ denote the set of functionals $g$ such that $|g(x + \delta y) - g(x)| \leq C\delta^\alpha$, for all $y \in L_2(\mathcal{I})$ satisfying $\|y\| \leq 1$, and for all $0 \leq \delta \leq 1$. When deriving convergence rates, we strengthen (9) by asking that $g$ be in $\mathcal{G}(C, x, \alpha)$.

Let $\mathcal{X} = \{X_1, \ldots, X_n\}$ denote the set of explanatory variables.

THEOREM 1. *If Assumptions 1 and 2 hold, then $\hat{g}(x) \to g(x)$ in mean square, conditional on $\mathcal{X}$, and*

$$\sup_{g \in \mathcal{G}(C,x,\alpha)} E[\{\hat{g}(x) - g(x)\}^2 | \mathcal{X}] = o_p(1). \tag{11}$$

*Furthermore, for all $\eta > 0$,*

$$\sup_{g \in \mathcal{G}(C,x,\alpha)} P\{|\hat{g}(x) - g(x)| > \eta\} \to 0.$$

*Moreover, if $h$ is chosen to decrease to zero in such a manner that*

$$h^{2\alpha} P(\|X - x\| \leq c_1 h) \asymp n^{-1} \tag{12}$$

*as $n \to \infty$, then, for each $C > 0$, the rate of convergence of $\hat{g}(x)$ to $g(x)$ equals $O_p(h^{2\alpha})$, uniformly in $g \in \mathcal{G}(C, x, \alpha)$:*

$$\sup_{g \in \mathcal{G}(C,x,\alpha)} E[\{\hat{g}(x) - g(x)\}^2 | \mathcal{X}] = O_p(h^{2\alpha}), \tag{13}$$

$$\lim_{C_1 \to \infty} \limsup_{n \to \infty} \sup_{g \in \mathcal{G}(C,x,\alpha)} P\{|\hat{g}(x) - g(x)| > C_1 h^\alpha\} = 0. \tag{14}$$

To interpret (11) and (13), assume that the pairs $(X_i, \varepsilon_i)$, for $1 \leq i < \infty$, are all defined on the same probability space, and then put $Y_i = Y_i(g) = g(X_i) + \varepsilon_i$. Write $E_g(\cdot | \mathcal{X})$ to denote expectation in the distribution of the pairs $(X_i, Y_i(g))$, conditional on $\mathcal{X}$. In Section 5.1 below, we shall discuss appropriateness of conditions such as (12), which relate to "small ball probabilities." Asymptotic consistency results for $g$ and mean squared errors have been derived in Ferraty, Mas and Vieu (2007) under different assumptions. The convergence rate at (14) is optimal in the following sense.

THEOREM 2. *If the error $\varepsilon$ in (1) is normally distributed, and if, for a constant $c_1 > 0$, $nP(\|X - x\| \leq c_1 h) \to \infty$ and (12) holds, then, for any estimator $\tilde{g}(x)$ of $g(x)$, and for $C > 0$ sufficiently large in the definition of $\mathcal{G}(C, x, \alpha)$, there exists a constant $C_1 > 0$, such that*

$$\limsup_{n \to \infty} \sup_{g \in \mathcal{G}(C,x,\alpha)} P\{|\tilde{g}(x) - g(x)| > C_1 h^\alpha\} > 0.$$

According to this result, uniformity of the convergence holds over the Lipschitz class of functionals $\mathcal{G}(C, x, \alpha)$. This result applies for a fixed argument $x$ in the domain of the predictor functions, where the functionals are evaluated. Further discussion of the bounds on $P(\|X - x\| \leq u)$ as relevant for (12) is provided in Section 5.1.

3.2. *Consistency of derivative estimator.* We shall establish consistency of the estimator $\hat\gamma_{xj}$. To this end, let

$$q_{12j} = 1 - \frac{|\int(X_1 - X_2)\psi_j|^2}{\|X_1 - X_2\|^2}$$

denote the version of $Q_{12j}$ when $\hat\xi_{i_1 i_2}$ is replaced by the quantity $\xi_j$ that $\hat\xi_{i_1 i_2}$ approximates, and let $k_{i_1 i_2 j}$ denote the version of $K(i_1, i_2, j|x)$, defined at (8), when $Q_{i_1 i_2 j}$ there is replaced by $q_{i_1 i_2 j}$.

ASSUMPTION 3.

(a) $\sup_{t \in \mathcal{I}} E\{X(t)^4\} < \infty$;
(b) there are no ties among the eigenvalues $\theta_1, \ldots, \theta_{j+1}$;
(c) $|g(x + y) - g(x) - g_x y| = o(\|y\|)$ as $\|y\| \to 0$;
(d) the distribution of $\xi_{1j} - \xi_{2j}$ has a well-defined density in a neighborhood of the origin, not vanishing at the origin;
(e) $K$ is supported on $[0, 1]$, nondecreasing and with a bounded derivative on the positive half-line, and $0 < K(0) < \infty$; and
(f) $h_1, h_2 \to 0$ as $n$ increases, sufficiently slowly to ensure that $n^{1/2} \min(h_1, h_2) \to \infty$ and $(nh_1)^2 E(k_{i_1 i_2 j}) \to \infty$.

Finite variance of $X$ guarantees that the covariance operator $V$, leading to the eigenfunctions $\psi_j$ and their estimators $\hat\psi_j$ in Section 3.1, is well defined; and finite fourth moment, stipulated by Assumption 4(a), ensures that $\|\hat\psi_j - \psi_j\|$ converges to zero at the standard root-$n$ rate. This assumption is, for example, satisfied for Gaussian processes with smooth mean and covariance functions.

If we suppose, in addition, that $X$ is a process with independent principal component scores $\int X\psi_j$ (or the stronger assumption that $X$ is Gaussian) and all the eigenvalues $\theta_j$ are nonzero [we shall refer to these properties jointly as $(P_1)$], then Assumption 3(f) implies that $n^{-\varepsilon} = O(h_j)$ for $j = 1, 2$ and for all $\varepsilon > 0$ [call this property $(P_2)$]. That is, both bandwidths are of larger order than any polynomial in $n^{-1}$. To see why, note that $(P_1)$ entails $P(\|x - X\| \leq h_1) = O(h_1^{C_1})$ for all $C_1 > 0$. Also, 3(f) implies that $nh_1 P(\|x - X\| \leq C_2 h_1) \to \infty$ for some $C_2 > 0$, and this, together with $(P_1)$, leads us to conclude that $nh^{C_1+1} \to \infty$ for all $C_1 > 0$. That result is equivalent to $(P_2)$ for the bandwidth $h_1$. Property $(P_1)$ also implies that $P(q_{12j} \leq h_2) = O(h_2^{C_1})$ for all $C_1 > 0$, and 3(f) implies that $nP(q_{12j} \leq C_2 h_2) \to \infty$ for some $C_2 > 0$, which, as before, leads to $(P_1)$, this time for the second bandwidth.

THEOREM 3. *If Assumption 3 holds, then $\hat\gamma_{xj} \to \gamma_{xj}$ in probability.*

Using notation (5), if $e = \sum_{j=1}^{j_0} e_j \psi_j$ with $\sum_j e_j^2 = 1$ and $j_0 < \infty$, the functional directional derivative in direction $e$ at $x$ is $g_x e = \sum_j e_j \gamma_{xj}$; see also (6), where $e$ is obtained by choosing $a_j = e_j, 1 \leq j \leq j_0, a_j = 0, j > j_0$. If Assumption 3 holds for all $j \leq j_0$, it is an immediate consequence of Theorem 3 that the estimated functional derivative $\hat{g}_x e = \sum_j e_j \hat{\gamma}_{xj}$ at $x$ in direction $e$ is consistent (i.e., satisfies $\hat{g}_x e \to g_x e$ in probability). As this holds uniformly over all direction vectors $e$, the functional gradient field for directions anchored in the span of $\{\psi_1, \ldots, \psi_{j_0}\}$ can be estimated consistently.

If we take the operator $\hat{g}_x$, defined by $\hat{g}_x a = \sum_{j \leq r} \hat{\gamma}_{xj} a_j$ (where $r \geq 1$ is an integer and $a = \sum_j a_j \psi_j$ is function), to be an empirical approximation to $g_x$, the operator given by $g_x a = \sum_j \gamma_{xj} a_j$, if the conditions in Assumption 3 hold for each $j$, and in addition $\sum_j \gamma_{xj}^2 < \infty$, then there exists a (generally unknown) deterministic sequence $r = r(n,x)$ with the following properties: $r(n,x) \to \infty$ as $n \to \infty$; whenever $\|a\| < \infty$, $\hat{g}_x a - g_x a \to 0$ in probability; and moreover, $\hat{g}_x \to g_x$ in norm as $n \to \infty$, where the convergence is again in probability. An explicit construction of such a sequence $r(n,x)$, and thus of an explicit estimate of the derivative operator with these properties, would require further results regarding the convergence rates for varying $j$ in Theorem 3, and remains an open problem.

**4. Application of functional derivative estimation to growth data.** The analysis of growth data has a long tradition in statistics. It played a pioneering role in the development of functional data analysis, as evidenced by the studies of Rao (1958), Gasser et al. (1984), Kneip and Gasser (1992), Ramsay and Li (1998) and Gervini and Gasser (2005) and remains an active field of statistical research to this day.

We explore the relationship between adult height, measured at age 18 (scalar response), and the growth rate function observed to age 10 (functional predictor), for 39 boys. Of interest is the following question: how do shape changes in the prepubertal growth velocity curve relate to changes in adult height? Which changes in the shape of a prepubertal growth velocity curve of an individual will lead to the largest adult height gain for an individual? These and similar questions can be addressed by obtaining the functional gradient of the regression of adult height versus the prepubertal growth velocity trajectory. Such analyses are expected to provide us with better understanding of the intricate dynamics and regulatory processes of human growth. Functional differentiation provides an excellent vehicle for studying the effects of localized growth velocity changes during various stages of prepubertal growth on adult height.

For this exploration, we use growth data for 39 boys from the Berkeley longitudinal growth study [Tuddenham and Snyder (1954)], where we include only measurements obtained up to age 10 for the growth velocity predictor processes. The 15 time points before age 10 at which height measurements are available for each boy in the Berkeley study correspond to ages $\{1, 1.25, 1.5, 1.75, 2, 3, 4, 5, 6, 7, 8, 8.5, 9, 9.5, 10\}$, denoted by $\{s_j\}_{j=1,\ldots,15}$. Raw growth rates were calculated as first order difference quotients $X_{ij} = (h_{i,j+1} - h_{ij})/(t_{j+1} - t_j)$, where $h_{ij}$ are the observed heights at times $s_j$ for the $i$th boy, and $t_j = (s_j + s_{j+1})/2$, $i = 1, \ldots, 39$, $j = 1, \ldots, 14$. These raw data form the input for the computation of the functional decomposition of the predictor processes into mean function, eigenfunctions and functional principal component scores. To obtain this decomposition, we used an implementation of the functional spectral methods described in Yao et al. (2003) and Yao, Müller and Wang (2005a).

Applying a BIC type criterion based on marginal pseudo-likelihood to choose the number of components in the eigenrepresentation, three components were selected. The resulting smooth estimates of fitted individual and mean growth velocity curves are shown in Figure 1. The first three components explain 99.5% of the total variation (78.9%, 17% and 3.6%, resp.), and the corresponding estimated eigenfunctions are displayed in the left panel of Figure 2. The first eigenfunction corresponds to a rapid initial decline in growth velocity, followed by a relatively flat increase with onset around age 5 toward the right end of the considered age range, while the second eigenfunction contains a sign change and provides a contrast between growth rates after age 2 and those before age 2. The third eigenfunction describes a midgrowth spurt around ages 6–7, coupled with an especially rapid decline in growth rate before age 3.

To visualize the estimated functional derivatives, a derivative scores plot as shown in the right panel of Figure 2 is of interest. The coefficient estimates for the first two eigendirections are plotted [i.e., the points, defined at (5), $(\gamma_{X_i,1}, \gamma_{X_i,2})$, evaluated at each of the 39 predictor functions $X_i$]. This figure thus represents the canonical functional gradient vectors at the observed data points, truncated at the first two components. These gradient vectors are seen to vary quite a bit across subjects, with a few extreme values present in the derivative corresponding to the first eigendirection.

The gradients are generally positive in the direction of the first eigenfunction and negative in the direction of the second. Their interpretation is relative to the shape of the eigenfunctions, including the selected sign for the eigenfunctions (as the sign of the eigenfunctions is arbitrary). If the gradient is positive in the direction of a particular eigenfunction $\psi_j$, it means that adult height tends to increase as the corresponding functional principal component score $\xi_j$ increases. So, in order to interpret the gradients in the right panel of Figure 2, one needs to study the shapes of the corresponding eigenfunctions as depicted in the left panel. When observing the shapes of first and second eigenfunction in the left panel of Figure 2, adult height is seen to increase most if the growth velocities toward the right end of the domain of the growth rate, predictor curves are larger, a result that is in line with what one would expect.

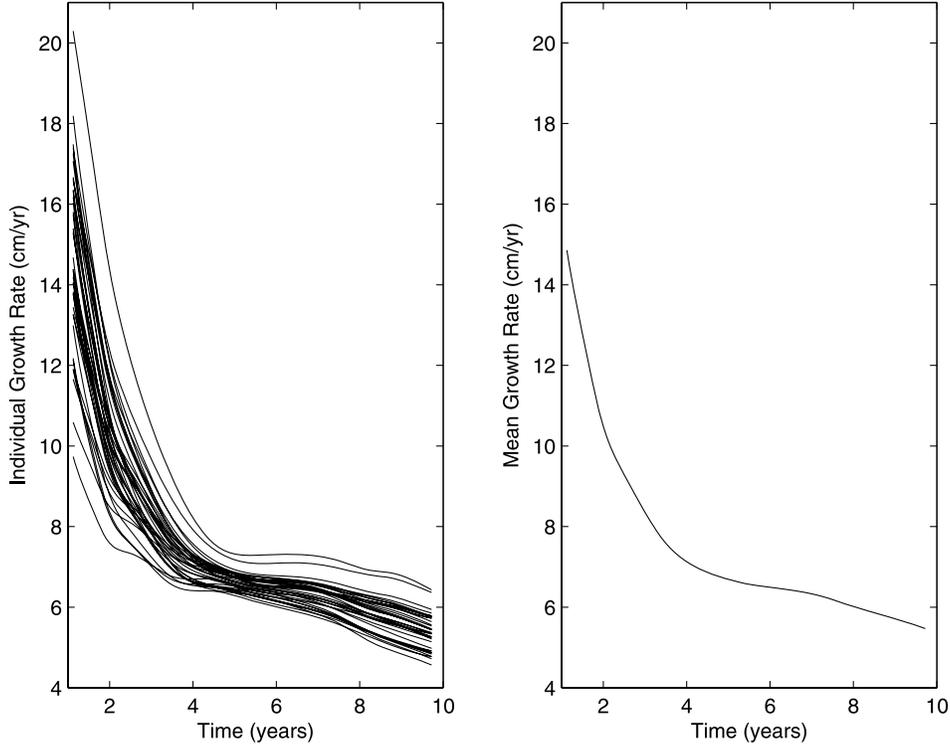

Fig. 1. *Fitted trajectories for individual predictor growth velocity curves (left panel) and mean growth velocity curve (right panel) for the Berkeley growth data ($n = 39$).*

Using the first $K$ components, we define functions $g_i^*(t) = \sum_{j=1}^{K} \gamma_{X_i,j} \psi_j(t)$ for each subject $i$. Then, for any test function $z(t) = \sum_{j=1}^{K} z_j \psi_j(t)$ with $\|z\| = 1$, one has $\int g_i^*(t) z(t)\, dt = \sum_{j=1}^{K} \gamma_{X_i,j} z_j$, so that the functional directional derivative at $X_i$ in direction $z$ is obtained through an inner product of $z$ with $g_i^*$. We therefore refer to $g_i^*$ as the *derivative generating function* at $X_i$. In the application to growth curves, we choose $K = 3$ and this function can be interpreted as a subject-specific weight function, whose inner product with a test function $z$ provides the gradient of adult height when moving from the trajectory $X_i$ in the direction indicated by $z$. It is straightforward to obtain estimates

$$\hat{g}_i^*(t) = \sum_{j=1}^{K} \hat{\gamma}_{X_i,j} \hat{\psi}_j(t) \tag{15}$$

of these derivative generating functions by plugging in estimates for $\gamma_{X_i,j}$ and $\psi_j(t)$ as obtained in (3) and (7).

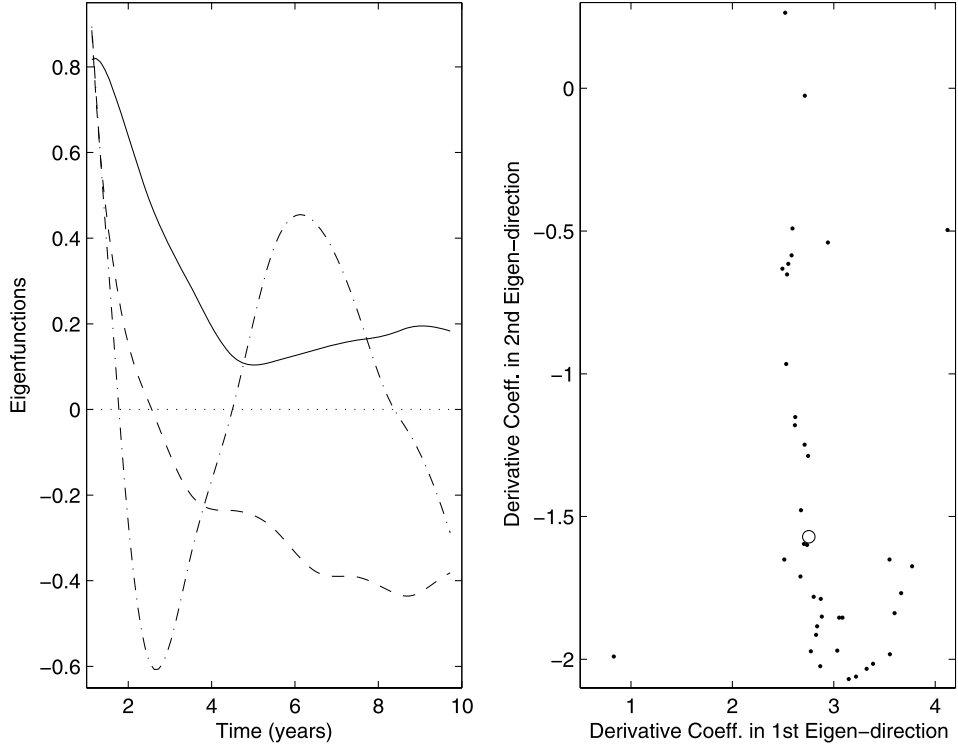

Fig. 2. *Smooth estimates of the first three eigenfunctions for the velocity growth curves, explaining 78.9% (solid), 17% (dashed) and 3.6% (dash-dotted) of the total variation, respectively (left panel) and estimated functional derivative coefficients $(\hat\gamma_{X_i,1}, \hat\gamma_{X_i,2})$ (7), in the directions of the first (x-axis) and second (y-axis) eigenfunction, evaluated at the predictor curves $X_i$ (dots), as well as at the mean curve $\mu$ (circle) (right panel).*

Estimated derivative generating functions $\hat g_i^*$ for $K = 3$ are depicted in Figure 3 for all 39 trajectories $X_i$ in the sample. These empirical derivative generating functions are found to be relatively homogeneous. Estimated functional directional derivatives in any specific direction of interest are then easily obtained. We find that gradients are largest in directions $z = g_i^*/\|g_i^*\|$ (i.e., in directions that are parallel to the derivative generating functions $g_i^*$). This means that largest increases in adult height are obtained in the presence of increased growth velocity around 2–4 years and past 8 years, while growth velocity increases between 5–7 years have only a relatively small effect.

It is of interest to associate the behavior of the derivative operators with features of the corresponding predictor trajectories. The predictor trajectories $X_i$, for which the derivative coefficients $\gamma_{X_i,j}$ have the largest and smallest absolute values in each of the first three eigendirections (for $j = 1, 2, 3$), are depicted in the upper panels of Figure 4. The lower panels show the corresponding derivative generating functions. One finds that the functional gradients of growth velocity curves that contain time periods of relatively small growth velocity are such that increased growth velocity in these time periods is associated with the largest increases in subsequent adult height (dashed curves in left and middle panel, dotted curve in right panel), as does slowing of above-normal high post-partum growth velocities (dashed curve in right panel).

A systematic visualization of the connection of predictor functions and the gradient field, as represented by the derivative generating functions, is obtained by considering families of predictor trajectories $X(t; \alpha_j) = \hat{\mu}(t) + \alpha_j \hat{\psi}_j(t)$ that move away from the mean growth velocity trajectory in the direction of a specific eigenfunction, while the other eigenfunctions are ignored, as shown in the upper panels of Figure 5 for the first three eigenfunctions. The corresponding derivative generating functions are in the lower panels. This visually confirms that adult height gains are associated with increased growth velocities in those areas where a subject's velocities are relatively

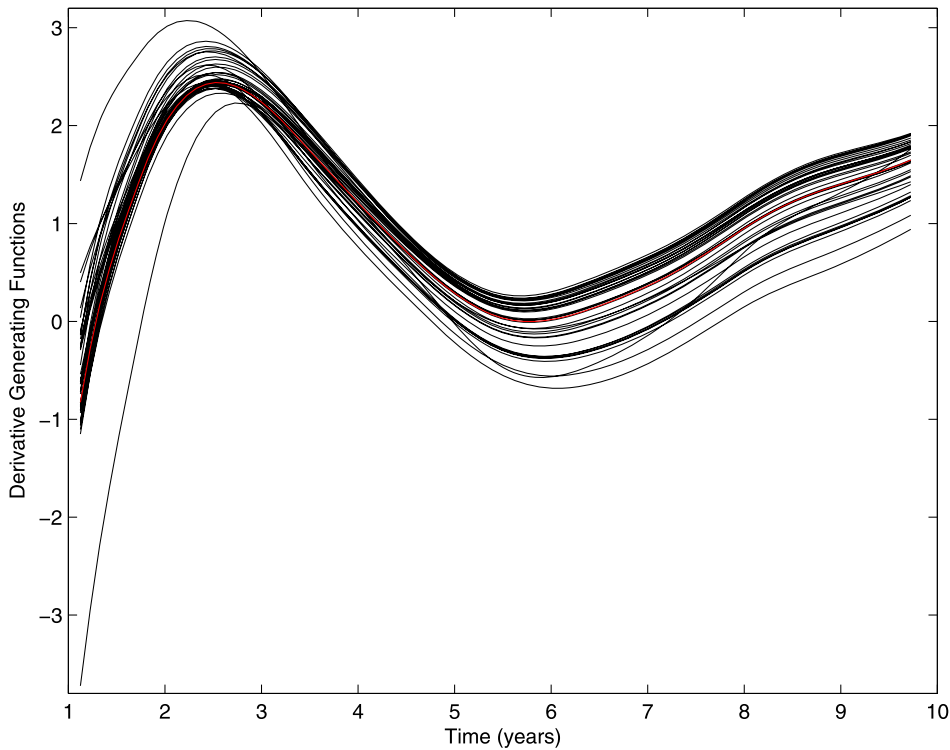

FIG. 3. *Estimated derivative generating functions $\hat{g}_i^*(t)$ (15) for all subjects $X_i$ (black) and for the mean function (red) of the Berkeley growth data, based on the first three eigenfunctions.*

low, especially toward the right end of the domain of the velocity predictor curves.

As the sample size in this example is relatively small, it is clear that caution needs to be exercised in the interpretation of the results of this data analysis. The results presented here follow the spirit of exploratory data analysis. We find that the concept of functional derivatives can lead to new insights when analyzing functional data which extend beyond those available when using established functional methods. Many practical and theoretical issues require further study. These include, for example, choice of window widths and the estimation of functional derivatives for data that are irregularly or sparsely measured.

## 5. Additional results and proofs.

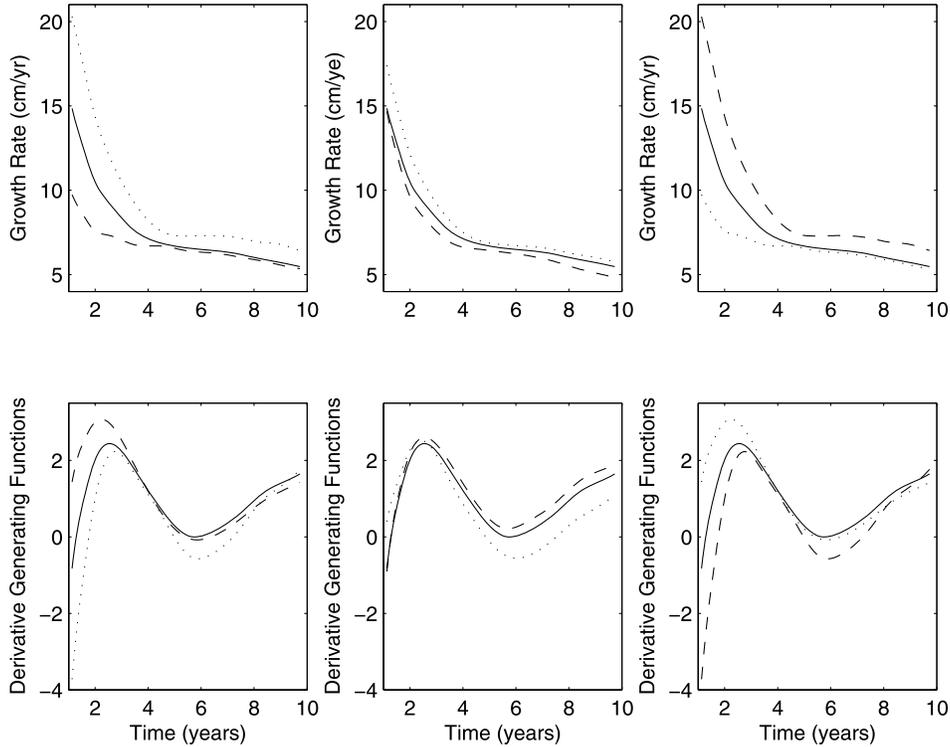

FIG. 4. *Predictor trajectories (top panels) and corresponding derivative generating functions $\hat{g}_i^*(t)$ (15) (bottom panels) which have the largest (dashed) and smallest (dotted) absolute values of derivative coefficients $\hat{\gamma}_{xj}$ (7) in the directions of the first ($j = 1$, left), second ($j = 2$, middle) and third ($j = 3$, right) eigenfunctions, as well as the mean functions (solid).*

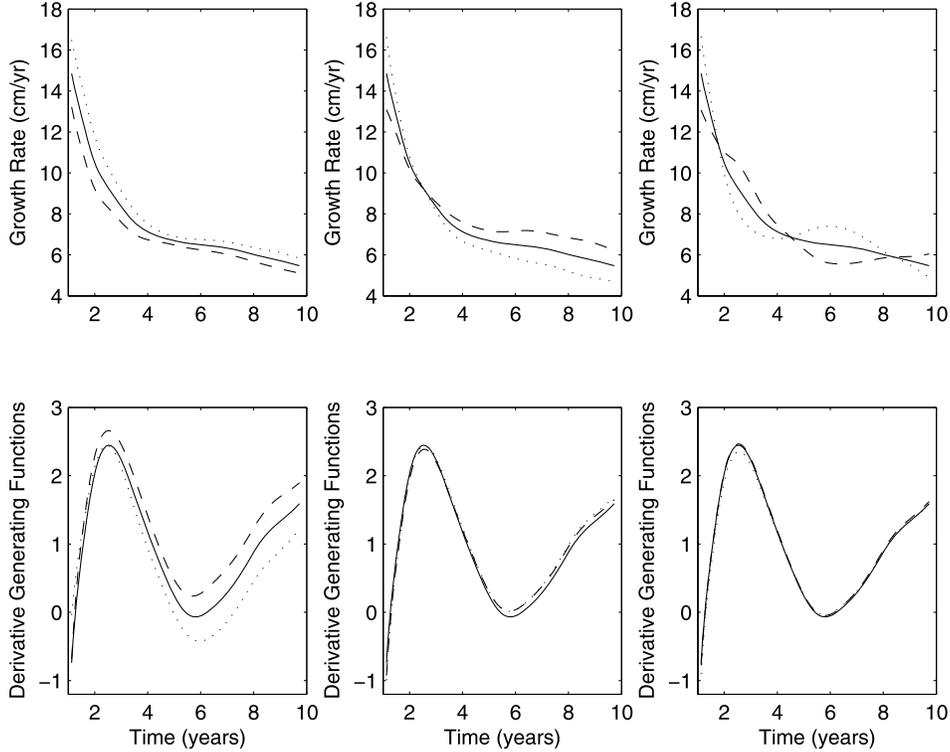

Fig. 5. *Top: predictor trajectories $X(t;\alpha_j) = \hat\mu(t) + \alpha_j \hat\psi_j(t)$ with $\alpha_j = -2$ (dashed), 0 (solid), +2 (dotted), where $j = 1,2,3$ from left to right. Bottom: corresponding derivative generating functions (15).*

5.1. *Bounds on $P(\|X - x\| \leq u)$.* Reflecting the infinite-dimensional nature of functional-data regression, the rate of convergence of the "small ball probabilities" $P(\|X - x\| \leq u)$ to zero as $u \to 0$ is generally quite rapid; in fact, it is faster than any polynomial in $u$. See (19) below. In consequence, the convergence rate of $\hat g(x)$ to $g(x)$ can be particularly slow. Indeed, unless the Karhunen–Loève expansion of $X$ is actually finite dimensional, the rate of convergence evidenced by (14) is slower than the inverse of any polynomial in $n$.

The fastest rates of convergence arise when the distribution of $X$ is closest to being finite dimensional; for example, when the Karhunen–Loève expansion of $X$ can be written as $X = \sum_j \xi_j \psi_j$, where $\operatorname{var}(\xi_j) = \theta_j$ and the eigenvalues $\theta_j, j \geq 1$, decrease to zero exponentially, rather than polynomially, fast as $j$ increases, where the $\xi_j$ are uncorrelated. Therefore, we shall focus primarily on this case and require the following.

ASSUMPTION 4. For constants $B, \beta > 0$,

$$\log \theta_j = -Bj^\beta + o(j^\beta) \qquad \text{as } j \to \infty, \tag{16}$$

and the random variables $\eta_j = \xi_j / \theta_j^{1/2}$ are independent and identically distributed as $\eta$, the distribution of which satisfies

$$B_1 u^b \leq P(|\eta| \leq u) \leq B_2 u^b \qquad \text{for all sufficiently small } u > 0, \quad \text{and}$$
$$P(|\eta| > u) \leq B_3 (1+u)^{-B_4} \qquad \text{for all } u > 0, \text{ where } B_1, \ldots, B_4, b > 0. \tag{17}$$

Take $x = 0$, the zero function, and, with $b$, $B$ and $\beta$ as in (16) and (17), define

$$\pi(u) = \exp\left\{-\frac{b\beta}{\beta+1}\left(\frac{2}{B}\right)^{1/\beta} |\log u|^{(\beta+1)/\beta}\right\}. \tag{18}$$

THEOREM 4. *If (16) and (17) hold, then, with $\pi(u)$ given by (18),*

$$P(\|X\| \leq u) = \pi(u)^{1+o(1)} \qquad \text{as } u \downarrow 0. \tag{19}$$

Combining Theorems 1 and 3, we deduce that, if the eigenvalues $\theta_j$ decrease as indicated at (16), if the principal components $\xi_j$ have the distributional properties at (17), and if the bandwidth $h$ is chosen so that (12) holds, then the kernel estimator $\hat{g}(x)$ converges to $g(x)$ at the mean-square rate of

$$h^{2\alpha} = \exp(-2\alpha|\log h|)$$
$$= \exp\left[-\{1+o(1)\}2\alpha\left(\frac{\beta+1}{b\beta}\right)^{\beta/(\beta+1)}\left(\frac{B}{2}\right)^{1/(\beta+1)} (\log n)^{\beta/(\beta+1)}\right].$$

For each fixed $\beta$, this quantity decreases to zero more slowly than any power of $n^{-1}$, although the rate of decrease increases as $\beta$ increases. A typical example where conditions (16) and (17) are satisfied is that of a process where $\theta_j = \exp(-Bj^\beta)$, where the distribution of $\eta$ in Assumption 4 has a bounded nonzero density in a neighborhood of the origin, and where $\{\phi_j\}$ is the standard Fourier series. In this case, one finds that $\beta = b = 1$ and $\pi(u) = \exp\{-c(\log u)^{(\beta+1)/\beta}\} = u^{-c(\log u)^{1/\beta}}$, for some $c > 0$, corresponding to faster than polynomial convergence toward 0. Of course, the condition on the distribution of $\eta$ is satisfied if the process $X$ is Gaussian.

Theorem 4 establishes that, in the case $x = 0$, the probability $P(\|X - x\| \leq u)$ typically does not vanish, even for very small $u$, and, in this context, (19) gives a concise account of the size of the probability. If we take $x = 0$ and replace $X$ by $X_1 - X_2$, for which the calculations leading to (19) are identical in all essential respects to those leading to (19), then we obtain a formula for the average value of $P(\|X_1 - x\| \leq u)$ over all realizations $x$ of $X_2$. Therefore,

(19) provides substantially more than just the value of the probability when $x = 0$. The case of fixed but nonzero $x$, where $x = \sum_j \theta_j^{1/2} x_j$ and the $x_j$'s are uniformly bounded, can be treated with related arguments, and also the setting where the $x_j$'s are unbounded, although it needs more detailed arguments.

If $\theta_j$ decreases to zero at a polynomial rate, rather than at the exponential rate stipulated by (16), then the probability $P(\|X - x\| \leq u)$ decreases to zero at rate $\exp(-C_1 u^{-C_2})$ as $u$ decreases to 0, rather than at the rate $\exp(-C_1 |\log u|^{C_2})$ indicated by Theorem 3 for constants $C_1, C_2 > 0$. Very accurate results of this type, in the case where $x = 0$, are given by Gao, Hannig and Torcaso (2003), who also provide additional relevant references. It is noteworthy that these results also pertain to non-Gaussian processes, while early results along these lines for Gaussian processes can be found in Anderson and Darling (1952). Decay rates of the closely related type $u^{C_3} \exp(-C_1 u^{-C_2})$ for $C_3 > 0$ were featured in Ferraty, Mas and Vieu (2007), among several other rates that are primarily associated with finite-dimensional processes.

We conclude from this discussion that the decay rates of the small ball probabilities are intrinsically linked to the decay rates of the eigenvalues of the underlying process. The fast decay rates associated with polynomially converging eigenvalues mean that this case is not particularly desirable from a statistical point of view.

5.2. *Proof of Theorem 1.* Let $\sigma^2$ denote the variance of the error $\varepsilon$ in (1). Set $N_j = \sum_i K_i(x)^j$ for $j = 1, 2$, and note that $N_2 \leq K(0) N_1$, as $K(\cdot)$ is nonincreasing and compactly supported on $[0, c]$. Therefore,

$$
\begin{aligned}
E[\{\hat{g}(x) - g(x)\}^2 | \mathcal{X}] \\
= [E\{\hat{g}(x) | \mathcal{X}\} - g(x)]^2 + \mathrm{var}(\hat{g}(x) | \mathcal{X}) \\
\leq \max_{i=1,\ldots,n} |g(X_i) - g(x)| I(\|X_i - x\| \leq ch) + \frac{\sigma^2 \sum_i K_i^2(x)}{\{\sum_i K_i(x)\}^2} \\
\leq \sup_{y:\,\|y\| \leq ch} |g(x) - g(x+y)|^2 + \frac{\sigma^2 K(0)}{N_1}.
\end{aligned}
\tag{20}
$$

Continuity of $g$ at $x$ [i.e., (9)] implies that the first term on the right-hand side of (20) converges to zero. Note that $K_i(x) \geq K_i(x) I(\|X_i - x\| \leq c_1 h) \geq K(c_1) I(\|X_i - x\| \leq c_1 h)$, where $c_1$ is as in (A2). Then, (10) entails $N_1^{-1} \to 0$ with probability 1, and by monotone convergence $E(N_i^{-1}) \to 0$. Together with (20), these properties imply the first part of the theorem. The second part, comprising (13) and (14), is obtained on noting that (20) entails

$$
\sup_{g \in \mathcal{G}(C,x,\alpha)} E[\{\hat{g}(x) - g(x)\}^2 | \mathcal{X}] \leq C^2 (ch)^{2\alpha} + \frac{\sigma^2 K(0)}{N_1}
$$

$$\leq C^2(ch)^{2\alpha} + \frac{\sigma^2 K(0)\{1+o_p(1)\}}{K(c_1)nP(\|X-x\|\leq c_1 h)}$$

and $E(N_1^{-1}) \leq E[\{\sum_i I(\|X_i - x\| \leq c_1 h)\}^{-1}] \asymp \{nP(\|X-x\|\leq c_1 h)\}^{-1}$.

5.3. *Proof of Theorem 2.* Without loss of generality, $x = 0$. Let $f$ denote a function defined on the real line, with a derivative bounded in absolute value by $B_1$, say, supported only within the interval $[-B_2, B_2]$, and not vanishing everywhere. Then, $f$ itself must be uniformly bounded, by $B_3$ say. Define $g_1 \equiv 0$ and $g_2(y) = h^\alpha f(\|y\|/h)$. If $\|y\| \leq h$ then, since $0 < \alpha \leq 1$,

$$|g_2(y) - g_2(0)| = h^\alpha |f(\|y\|/h) - f(0)| \leq h^\alpha B_1 \|y\|/h \leq h^\alpha B_1 (\|y\|/h)^\alpha$$
$$= B_1 \|y\|^\alpha,$$

while, if $\|y\| > h$,

$$|g_2(y) - g_2(0)| \leq 2h^\alpha B_3 \leq 2B_3 \|y\|^\alpha.$$

Therefore, $g_2 \in \mathcal{G}(C, 0, \alpha)$ provided $\max(B_1, 2B_3) \leq C$.

The theorem will follow if we show that, in a classification problem where we observe $n$ data generated as at (1), with the errors distributed as Normal $N(0, 1)$ and $g = g_1$ or $g_2$, with prior probability $\frac{1}{2}$ on either of these choices, the likelihood-ratio rule fails, in the limit as $n \to \infty$, to discriminate between $g_1$ and $g_2$. That is, with $Y_i = \varepsilon_i$ (the result of taking $g = g_1$ in the model), and with $\rho$ defined by

$$\rho = \frac{\prod_i \exp[-1/2\{Y_i - g_1(X_i)\}^2]}{\prod_i \exp[-1/2\{Y_i - g_2(X_i)\}^2]},$$

we should show that

(21) $\qquad P(\rho > 1)$ is bounded below 1 $\qquad$ as $n \to \infty$.

Now,

$$2 \log \rho = \sum_{i=1}^n \{g_2(X_i)^2 - 2\varepsilon_i g_2(X_i)\},$$

which, conditional on $\mathcal{X}$, is normally distributed with mean $s_n^2 = \sum_i g_2(X_i)^2$ and variance $4s_n^2$. Therefore, (21) holds if and only if

(22) $$\lim_{B \to \infty} \limsup_{n \to \infty} P(s_n^2 > B) = 0$$

and so we can complete the proof of Theorem 2 by deriving (22).

If we choose the radius $B_2$ of the support of $f$ so that $0 < B \leq c_1$, then $|g_2(x)| \leq B_3 h^\alpha I(\|x\| \leq c_1 h)$, in which case

(23) $$s_n^2 \leq B_3^2 h^{2\alpha} \sum_{i=1}^n I(\|X_i\| \leq c_1 h).$$

Since, by assumption, $nP(\|X\| \leq c_1 h) \to \infty$, then

$$\frac{\sum_i I(\|X_i\| \leq c_1 h)}{nP(\|X\| \leq c_1 h)} \to 1$$

in probability. This property, (12) and (23) together imply (22).

5.4. *Proof of Theorem 3.* Write, simply, $K_{i_1 i_2 j}$ for $K(i_1, i_2, j|x)$. Assumption 3(e) implies that

(24) $$K_{i_1 i_2 j} = 0, \quad \text{unless each of the following holds: } \|X_{i_1} - x\| \leq h_1,$$
$$\|X_{i_2} - x\| \leq h_1 \text{ and } Q_{i_1 i_2} \leq h_2.$$

Given $\delta > 0$, let $s(\delta)$ equal the supremum of $|g(x+y) - g(x) - g_x y|$ over functions $y$ with $\|y\| \leq \delta$. Then, by Assumption 3(c),

(25) $$\delta^{-1} s(\delta) \to 0 \quad \text{as } \delta \downarrow 0.$$

Write $\mathcal{E}_{i_1 i_2}$ for the event that $\|X_{i_k} - x\| \leq h_1$ for $k = 1, 2$. If $\mathcal{E}_{i_1 i_2}$ holds,

$$|g(X_{i_1}) - g(X_{i_2}) - g_x(X_{i_1} - X_{i_2})| \leq 2s(h_1).$$

Therefore, defining $\varepsilon_{i_1 i_2} = \varepsilon_{i_1} - \varepsilon_{i_2}$ and assuming $\mathcal{E}_{i_1 i_2}$,

$$|Y_{i_1} - Y_{i_2} - \{g_x(X_{i_1} - X_{i_2}) + \varepsilon_{i_1 i_2}\}| \leq 2s(h_1).$$

Hence, defining $\xi_{i_1 i_2 j} = \xi_{i_1 j} - \xi_{i_2 j}$, noting that $g_x(X_{i_1} - X_{i_2}) = \sum_k \xi_{i_1 i_2 k} \gamma_{xk}$, and using (24), we have,

(26) $$\left| \sum_{i_1, i_2}^{(j)} (Y_{i_1} - Y_{i_2}) K_{i_1 i_2 j} \right.$$
$$\left. - \left( \sum_{i_1, i_2}^{(j)} K_{i_1 i_2 j} \sum_{k=1}^{\infty} \xi_{i_1 i_2 k} \gamma_{xk} + \sum_{i_1, i_2}^{(j)} \varepsilon_{i_1 i_2} K_{i_1 i_2 j} \right) \right|$$
$$\leq 2s(h_1) \sum_{i_1, i_2}^{(j)} K_{i_1 i_2 j}.$$

Now,

(27) $$|\hat{\xi}_{i_1 i_2 j} - \xi_{i_1 i_2 j}| = \left| \int (X_{i_1} - X_{i_2})(\hat{\psi}_j - \psi_j) \right|$$
$$\leq \|X_{i_1} - X_{i_2}\| \|\hat{\psi}_j - \psi_j\| \leq 2h_1 \|\hat{\psi}_j - \psi_j\|,$$

where the last inequality holds under the assumption that the event $\mathcal{E}_{i_1 i_2}$ obtains. Combining (24), (26) and (27), we deduce that

$$
(28) \quad \left| \sum_{i_1,i_2}^{(j)} (Y_{i_1} - Y_{i_2}) K_{i_1 i_2 j} - \left( \gamma_{xj} \sum_{i_1,i_2}^{(j)} \hat{\xi}_{i_1 i_2 j} K_{i_1 i_2 j} + \sum_{i_1,i_2}^{(j)} K_{i_1 i_2 j} \sum_{k : k \neq j} \xi_{i_1 i_2 k} \gamma_{xk} + \sum_{i_1,i_2}^{(j)} \varepsilon_{i_1 i_2} K_{i_1 i_2 j} \right) \right|
$$

$$
\leq 2 \{ s(h_1) + |\gamma_{xj}| h_1 \|\hat{\psi}_j - \psi_j\| \} \sum_{i_1,i_2}^{(j)} K_{i_1 i_2 j}.
$$

Note, too, that

$$
\left| \sum_{i_1,i_2}^{(j)} K_{i_1 i_2 j} \sum_{k : k \neq j} \xi_{i_1 i_2 k} \gamma_{xk} \right|
$$

$$
= \left| \sum_{i_1,i_2}^{(j)} K_{i_1 i_2 j} \sum_{k : k \neq j} \gamma_{xk} \int (X_{i_1} - X_{i_2}) \psi_k \right|
$$

$$
(29) \quad \leq \sum_{i_1,i_2}^{(j)} K_{i_1 i_2 j} \left( \sum_{k : k \neq j} \gamma_{xk}^2 \right)^{1/2} \left[ \sum_{k : k \neq j} \left\{ \int (X_{i_1} - X_{i_2}) \psi_k \right\}^2 \right]^{1/2}
$$

$$
\leq \|g_x\| \sum_{i_1,i_2}^{(j)} K_{i_1 i_2 j} \left[ \|X_{i_1} - X_{i_2}\|^2 - \left\{ \int (X_{i_1} - X_{i_2}) \psi_j \right\}^2 \right]^{1/2}
$$

$$
\leq \|g_x\| \sum_{i_1,i_2}^{(j)} K_{i_1 i_2 j} \left[ \|X_{i_1} - X_{i_2}\|^2 - \left\{ \int (X_{i_1} - X_{i_2}) \hat{\psi}_j \right\}^2 \right.
$$

$$
\left. + 8 \|\hat{\psi}_j - \psi_j\| \|X_{i_1} - X_{i_2}\|^2 \right]^{1/2}
$$

$$
\leq 2 \|g_x\| h_1 \sum_{i_1,i_2}^{(j)} K_{i_1 i_2 j} (Q_{i_1 i_2 j} + 8 \|\hat{\psi}_j - \psi_j\|)^{1/2}
$$

$$
\leq 2 \|g_x\| h_1 (h_2 + 8 \|\hat{\psi}_j - \psi_j\|)^{1/2} \sum_{i_1,i_2}^{(j)} K_{i_1 i_2 j}.
$$

To obtain the third-last inequality in (29), we used the fact that, with $a = |\int (X_{i_1} - X_{i_2})\psi_j|$, $b = |\int (X_{i_1} - X_{i_2})\hat{\psi}_j|$ and

$$(30) \qquad c = \|X_{i_1} - X_{i_2}\|\|\hat{\psi}_j - \psi_j\| \leq 2\|X_{i_1} - X_{i_2}\| \leq 4h_1,$$

where [in each of (30) and in (31) below] the last inequality is correct provided $\mathcal{E}_{i_1 i_2}$ holds, we have used the fact that $|a-b| \leq c$ and $|a| \leq \|X_{i_1} - X_{i_2}\|$ imply that

$$(31) \qquad \begin{aligned} |a^2 - b^2| &\leq c(2a + c) \leq 4\|\hat{\psi}_j - \psi_j\|\|X_{i_1} - X_{i_2}\|^2 \\ &\leq 8\|\hat{\psi}_j - \psi_j\|h_1^2. \end{aligned}$$

To obtain the last inequality in (29), we used (24) and the fact that $Q_{i_1 i_2 j} \leq h_2$ if $K_{i_1 i_2 j} \neq 0$.

Combining (28) and (29), we find that

$$(32) \qquad \begin{aligned} &\left| \sum \sum_{i_1, i_2}^{(j)} (Y_{i_1} - Y_{i_2}) K_{i_1 i_2 j} \right. \\ &\qquad \left. - \left( \gamma_{xj} \sum \sum_{i_1, i_2}^{(j)} \hat{\xi}_{i_1 i_2 j} K_{i_1 i_2 j} + \sum \sum_{i_1, i_2}^{(j)} \varepsilon_{i_1 i_2} K_{i_1 i_2 j} \right) \right| \\ &\quad \leq 2h_1 \{ h_1^{-1} s(h_1) + |\gamma_{xj}| \|\hat{\psi}_j - \psi_j\| \\ &\qquad + \|g_x\|(h_2 + 8\|\hat{\psi}_j - \psi_j\|)^{1/2} \} \sum \sum_{i_1, i_2}^{(j)} K_{i_1 i_2 j}. \end{aligned}$$

Result (32) controls the numerator in the definition of $\hat{\gamma}_{xj}$ at (7). To control the denominator there, use (27) to show that

$$(33) \qquad \begin{aligned} \sum \sum_{i_1, i_2}^{(j)} \hat{\xi}_{i_1 i_2 j} K_{i_1 i_2 j} &\geq \sum \sum_{i_1, i_2}^{(j)} \max(0, \xi_{i_1 j} - \xi_{i_2 j} - 2h_1 \|\hat{\psi}_j - \psi_j\|) K_{i_1 i_2 j} \\ &\geq \sum \sum_{i_1, i_2}^{(j)} \max(0, \xi_{i_1 j} - \xi_{i_2 j}) K_{i_1 i_2 j} \\ &\quad - 2h_1 \|\hat{\psi}_j - \psi_j\| \sum \sum_{i_1, i_2}^{(j)} K_{i_1 i_2 j}. \end{aligned}$$

[Recall that $\sum \sum_{i_1, i_2}^{(j)}$ denotes summation over $(i_1, i_2)$ such that $\hat{\xi}_{i_1 i_2 j} > 0$.] Using Assumption 4(d), (e) and (f), it can be proved that, for a constant

ESTIMATION OF FUNCTIONAL DERIVATIVES 21$B > 0$,

$$(34) \quad \sum_{i_1,i_2}^{(j)} \max(0, \xi_{i_1 j} - \xi_{i_2 j}) K_{i_1 i_2 j} \geq \{1 + o_p(1)\} B h_1 \sum_{i_1,i_2}^{(j)} K_{i_1 i_2 j}.$$

[Note that, by Assumption 3(f), $n^{-1/2}/\min(h_1, h_2) \to 0$.] From Assumption 3(a) and (b), it follows that

$$(35) \quad \|\hat{\psi}_j - \psi_j\| = O_p(n^{-1/2}).$$

Together, (33)–(35) imply that

$$(36) \quad \sum_{i_1,i_2}^{(j)} \hat{\xi}_{i_1 i_2 j} K_{i_1 i_2 j} \geq \{1 + o_p(1)\} B h_1 \sum_{i_1,i_2}^{(j)} K_{i_1 i_2 j}$$

for the same constant $B$ as in (34). This result controls the denominator at (7).

From (7), (25), (32) and (36), we deduce that

$$(37) \quad \hat{\gamma}_{xj} = \gamma_{xj} + O_p\left(\frac{\sum \sum_{i_1,i_2}^{(j)} \varepsilon_{i_1 i_2} K_{i_1 i_2 j}}{h_1 \sum \sum_{i_1,i_2}^{(j)} K_{i_1 i_2 j}}\right) + o_p(1).$$

The variance of the ratio on the right-hand side of (37), conditional on the explanatory variables $X_i$, equals

$$O_p\left\{\left(h_1^2 \sum_{i_1,i_2}^{(j)} K_{i_1 i_2 j}\right)^{-1}\right\} = O_p[\{(nh_1)^2 E(k_{i_1 i_2 j})\}^{-1}] = o_p(1),$$

where, to obtain the last identity, we used Assumption 3(f). Therefore, (37) implies that $\hat{\gamma}_{xj} = \gamma_{xj} + o_p(1)$, which proves Theorem 3.

5.5. *Proof of Theorem 4.* Observe that, for each $t \in (0,1)$ and with $D_t = (\sum_j \theta_j^{1-t})^{-1}$,

$$(38) \quad P(\|X\| \leq u) = P\left(\sum_{j=1}^{\infty} \theta_j \eta_j^2 \leq u^2\right) \begin{cases} \leq \prod_{j=1}^{\infty} P(\theta_j \eta_j^2 \leq u^2), \\ \geq \prod_{j=1}^{\infty} P(\theta_j^t \eta_j^2 \leq D_t u^2), \end{cases}$$

where, to obtain the lower bound, we used the property

$$P\left(\sum_{j=1}^{\infty} \theta_j \eta_j^2 \leq u^2\right) = P\left\{\sum_{j=1}^{\infty} \theta_j^{1-t}(\theta_j^t \eta_j^2 - D_t u^2) \leq 0\right\}$$

$$\geq P(\theta_j^t \eta_j^2 \leq D_t u^2 \text{ for each } j).$$

Define $J = J(u)$ to be the largest integer such that $u/\theta_j^{1/2} \leq \zeta$, where $\zeta$ is chosen so small that $B_1 u^b \leq P(|\eta| \leq u) \leq B_2 u^b$ for $0 \leq u \leq \zeta$. Then,

$$
\begin{aligned}
\prod_{j=1}^{\infty} P(\theta_j \eta_j^2 \leq u^2) &\leq \prod_{j=1}^{J} P(|\eta| \leq u\theta_j^{-1/2}) \\
&= u^{bJ} \exp\left\{\frac{1}{2}bB \sum_{j=1}^{J} j^{\beta} + o(J^{\beta+1})\right\} \\
&= \exp\left\{-\frac{bB\beta}{2(\beta+1)} J^{\beta+1} + o(J^{\beta+1})\right\} \\
&= \pi(u)^{1+o(1)}
\end{aligned}
\tag{39}
$$

as $u \downarrow 0$, where $\pi$ is defined at (18).

Redefine $J$ to be the largest integer such that $D_t^{1/2} u/\theta_j^{t/2} \leq \zeta$. Then, using the argument leading to (39), we may show that

$$
\begin{aligned}
\prod_{j=1}^{J} P(\theta_j^t \eta_j^2 \leq D_t u^2) &= \exp\left\{-\frac{b\beta}{\beta+1}\left(\frac{2}{Bt}\right)^{1/\beta} |\log u|^{(\beta+1)/\beta} + o(|\log u|^{(\beta+1)/\beta})\right\} \\
&= \pi(u)^{t^{-1/\beta}+o(1)}.
\end{aligned}
\tag{40}
$$

Also, for $j \geq J+1$,

$$
\pi_j \equiv P(\theta_j^t \eta_j^2 > D_t u^2) \leq B_3 \{1 + (D_t^{1/2} u/\theta_j^{t/2})\}^{-B_4}. \tag{41}
$$

Note, too, that, for a constant $B_5 = B_5(t) \in (0,1)$, we have $\pi_j \in (0, B_5)$ for $j \geq J+1$, and

$$
1 - \pi_j = \exp\left(-\sum_{k=1}^{\infty} \frac{\pi_j^k}{k}\right) \geq \exp(-B_6 \pi_j)
$$

from which it follows that

$$
\prod_{j=J+1}^{\infty} (1 - \pi_j) \geq \exp\left(-B_6 \sum_{j=J+1}^{\infty} \pi_j\right) \geq \exp\left\{-B_7 \sum_{j=J+1}^{\infty} (\theta_j^{t/2}/u)^{B_4}\right\},
$$

which is of smaller order than the right-hand side of (40). Combining this result with (40), and noting that $t \in (0,1)$ on the right-hand side of (40), can be taken arbitrarily close to 1, we deduce that, as $u \downarrow 0$,

$$
\prod_{j=1}^{\infty} P(\theta_j^t \eta_j^2 \leq D_t u^2) = \pi(u)^{1+o(1)}. \tag{42}
$$

Together, (38), (39) and (42) imply (19).

ESTIMATION OF FUNCTIONAL DERIVATIVES 23

**Acknowledgments.** We wish to thank an Associate Editor and two referees for helpful comments.

## REFERENCES

Anderson, T. and Darling, D. (1952). Asymptotic theory of certain "goodness of fit" criteria based on stochastic processes. *Ann. Math. Statistics* **23** 193–212. MR0050238

Cardot, H., Ferraty, F., Mas, A. and Sarda, P. (2003). Testing hypotheses in the functional linear model. *Scand. J. Statist.* **30** 241–255. MR1965105

Cardot, H., Ferraty, F. and Sarda, P. (2003). Spline estimators for the functional linear model. *Statist. Sinica* **13** 571–591. MR1997162

Cuevas, A., Febrero, M. and Fraiman, R. (2002). Linear functional regression: The case of fixed design and functional response. *Canad. J. Statist.* **30** 285–300. MR1926066

Ferraty, F., Mas, A. and Vieu, P. (2007). Nonparametric regression on functional data: Inference and practical aspects. *Aust. N. Z. J. Stat.* **49** 459–461. MR2396496

Ferraty, F. and Vieu, P. (2003). Curves discrimination: A nonparametric functional approach. *Comput. Statist. Data Anal.* **44** 161–173. MR2020144

Ferraty, F. and Vieu, P. (2004). Nonparametric models for functional data, with application in regression, time-series prediction and curve discrimination. *J. Nonparametr. Stat.* **16** 111–125. MR2053065

Ferraty, F. and Vieu, P. (2006). *Nonparametric Functional Data Analysis.* Springer, New York. MR2229687

Gasser, T., Hall, P. and Presnell, B. (1998). Nonparametric estimation of the mode of a distribution of random curves. *J. R. Stat. Soc. Ser. B Stat. Methodol.* **60** 681–691. MR1649539

Gasser, T., Köhler, W., Müller, H.-G., Kneip, A., Largo, R., Molinari, L. and Prader, A. (1984). Velocity and acceleration of height growth using kernel estimation. *Annals of Human Biology* **11** 397–411.

Gao, F., Hannig, J. and Torcaso, F. (2003). Comparison theorems for small deviations of random series. *Electron. J. Probab.* **8** 17pp. MR2041822

Gervini, D. and Gasser, T. (2005). Nonparametric maximum likelihood estimation of the structural mean of a sample of curves. *Biometrika* **92** 801–820. MR2234187

Hall, P. and Heckman, N. E. (2002). Estimating and depicting the structure of a distribution of random functions. *Biometrika* **89** 145–158. MR1888371

Hall, P. and Horowitz, J. L. (2007). Methodology and convergence rates for functional linear regression. *Ann. Statist.* **35** 70–91. MR2332269

James, G. and Silverman, B. (2005). Functional adaptive model estimation. *J. Amer. Statist. Assoc.* **100** 565–576. MR2160560

Kneip, A. and Gasser, T. (1992). Statistical tools to analyze data representing a sample of curves. *Ann. Statist.* **20** 1266–1305. MR1186250

Ramsay, J. O. and Dalzell, C. J. (1991). Some tools for functional data analysis (with discussion). *J. Roy. Statist. Soc. Ser. B* **53** 539–572. MR1125714

Ramsay, J. O. and Li, X. (1998). Curve registration. *J. R. Stat. Soc. Ser. B Stat. Methodol.* **60** 351–363. MR1616045

Ramsay, J. O. and Silverman, B. W. (2005). *Functional Data Analysis*, 2nd ed. Springer, New York. MR2168993

Rao, C. R. (1958). Some statistical methods for comparison of growth curves. *Biometrics* **14** 1–17.

Tuddenham, R. and Snyder, M. (1954). Physical growth of California boys and girls from birth to age 18. *Calif. Publ. Child Develop.* **1** 183–364.

YAO, F., MÜLLER, H.-G., CLIFFORD, A. J., DUEKER, S. R., FOLLETT, J., LIN, Y., BUCHHOLZ, B. A. and VOGEL, J. S. (2003). Shrinkage estimation for functional principal component scores with application to the population kinetics of plasma folate. *Biometrics* **59** 676–685. MR2004273

YAO, F., MÜLLER, H.-G. and WANG, J.-L. (2005a). Functional data analysis for sparse longitudinal data. *J. Amer. Statist. Assoc.* **100** 577–590. MR2160561

YAO, F., MÜLLER, H.-G. and WANG, J.-L. (2005b). Functional linear regression analysis for longitudinal data. *Ann. Statist.* **33** 2873–2903. MR2253106

P. HALL  
DEPARTMENT OF MATHEMATICS  
AND STATISTICS  
UNIVERSITY OF MELBOURNE  
PARKVILLE, VIC, 3010  
AUSTRALIA  

H.-G. MÜLLER  
DEPARTMENT OF STATISTICS  
UNIVERSITY OF CALIFORNIA  
ONE SHIELDS AVENUE  
DAVIS, CALIFORNIA 95616  
USA  
E-MAIL: mueller@wald.ucdavis.edu  

F. YAO  
DEPARTMENT OF STATISTICS  
UNIVERSITY OF TORONTO  
100 SAINT GEORGE STREET  
TORONTO, ONTARIO M5S3G3  
CANADA